\documentclass[12pt]{amsart}
\usepackage{amsfonts,latexsym,amsthm,amssymb,graphicx}
\usepackage[all]{xy}

\newtheorem{theorem}{Theorem}[section]
\newtheorem{lemma}[theorem]{Lemma}
\newtheorem{corol}[theorem]{Corollary}
\newtheorem{prop}[theorem]{Proposition}
\newtheorem{claim}[theorem]{Claim}

{\theoremstyle{definition} \newtheorem{defin}[theorem]{Definition}}
{\theoremstyle{remark} \newtheorem{remark}[theorem]{Remark}
}

\newcommand{\Abb}{{\mathbb{A}}}
\newcommand{\Gbb}{{\mathbb{G}}}
\newcommand{\Pbb}{{\mathbb{P}}}
\newcommand{\Qbb}{{\mathbb{Q}}}
\newcommand{\Zbb}{{\mathbb{Z}}}
\newcommand{\Til}[1]{{\widetilde{#1}}}

\DeclareMathOperator{\rk}{rk}
\DeclareMathOperator{\Spec}{Spec}
\DeclareMathOperator{\Grass}{Grass}

\newcommand{\one}{1\hskip-3.5pt1}
\newcommand{\bxd}[1]{\{{#1}\}}

\title[pro-Chow groups and pro-CSM classes]{Limits of Chow groups, 
and a new construction of 
Chern-Schwartz-MacPherson classes}
\author{Paolo Aluffi}

\dedicatory
{Dedicated to Robert MacPherson on the occasion of his 60th
birthday}

\address{Max-Planck-Institut f\"ur Mathematik, Bonn, Deutschland}
\address{Math Dept, Florida State University, Tallahassee, Florida, U.S.A.}
\email{aluffi@math.fsu.edu}

\begin{document}

\begin{abstract}
We define an `enriched' notion of Chow groups for algebraic varieties, agreeing
with the conventional notion for complete varieties, but enjoying a functorial
push-forward for arbitrary maps. This tool allows us to glue intersection-theoretic
information across elements of a stratification of a variety; we illustrate this 
operation by giving a direct construction of {\em Chern-Schwartz-MacPherson\/} 
classes of singular varieties, providing a new proof of an old (and long since
settled) conjecture of Deligne and Grothendieck.
\end{abstract}

\maketitle

\setcounter{tocdepth}{1}
\tableofcontents


\newcommand{\hA}{{\widehat A}}
\newcommand{\cC}{{\mathcal C}}
\newcommand{\fc}{{\mathfrak c}}
\newcommand{\oU}{{\overline U}}
\newcommand{\oV}{{\overline V}}
\newcommand{\vvv}{{V}}
\newcommand{\cF}{{\mathcal F}}
\newcommand{\cQ}{{\mathcal Q}}
\newcommand{\tgamma}{{\widetilde\gamma}}
\newcommand{\Eus}{{E_{\underline s}}}
\newcommand{\cSg}{{{\mathcal S}_\Gamma}}
\newcommand{\cQg}{{{\mathcal Q}_\Gamma}}
\newcommand{\cT}{{\mathcal T}}
\newcommand{\hcT}{{\widehat \cT}}
\newcommand{\tGamma}{{\Til\Gamma}}
\newcommand{\cK}{{\mathcal K}}

\section{Introduction}\label{intro}
In the remarkable article \cite{MR50:13587}, Robert MacPherson
settled affirmatively a conjecture of Pierre Deligne and Alexandre 
Grothendieck  (see \cite{MR0278333} p.~168; and \cite{ReS}, 
note $({}^{87_1})$, for Grothendieck's own comments on the genesis
of the original conjecture). MacPherson's theorem states that there is a 
unique natural transformation from a functor of constructible functions 
on compact complex algebraic varieties to homology, associating 
to the constant function~$\one_V$ on a nonsingular variety $V$ the
(Poincar\'e dual of the) total Chern class of the tangent bundle $TV$
of $V$. The class corresponding to the constant $\one_X$ for an {\em
arbitrary\/} compact complex algebraic variety $X$ is therefore a very 
natural candidate for a notion of {\em Chern class\/} of a possibly 
singular variety.

After MacPherson's work it was realized that these classes agree, up to 
Alexander duality, with classes defined earlier by Marie-H\'el\`ene 
Schwartz (\cite{MR35:3707}, \cite{MR32:1727}; and 
\cite{MR83h:32011}). It is common nowadays to name these classes 
{\em Chern-Schwartz-MacPherson\/} (CSM) classes. 

MacPherson's construction may be used to lift the classes to the Chow
group $A_*X$, cf.~\cite{MR85k:14004}, \S19.1.7. Several other
approaches to CSM classes are known: for example through local polar 
varieties (\cite{MR634426}); characteristic cycles and index formulas for
holonomic $\mathcal D$-modules (\cite{MR647684}, \cite{MR804052}, 
\cite{MR833194}); currents and curvature measures (\cite{MR1267891}). 
Some of these approaches may be used to extend the definition 
of CSM classes to varieties over arbitrary algebraically closed fields of 
characteristic~zero (see \cite{MR91h:14010}), proving naturality
at the level of Chow groups, under proper push-forward.

The main goal of this paper is to provide a new construction of CSM
classes in this algebro-geometric setting, independent of previous 
approaches, and including a complete
proof of the naturality mandated by the Deligne--Grothendieck
conjecture. Our approach is very direct: we simply define 
an invariant for nonsingular (but possibly noncomplete) 
varieties, and obtain the CSM class of an arbitrary variety~$X$ 
as the sum of these invariants, over any decomposition
of $X$ as finite disjoint union of nonsingular subvarieties.
The contribution for a nonsingular variety is obtained
as Chern class of the dual of a bundle of differential forms
with logarithmic poles.

The approach is particularly transparent, as the contribution of a
nonsingular subvariety $U$ to the class for $X$ is independent of how
singular $X$ is along $U$. Auxiliary invariants, such as the 
{\em local Euler obstruction\/} or the {\em Chern-Mather\/} class, 
which are common to several of the approaches listed above, play 
no r\^ole in our construction.

Naturality is straightforward, modulo one technical lemma 
(Lemma~\ref{keycov}; see also Claim~\ref{keyclaim}) on the behavior 
of the contributions under push-forward. The classes we define must 
agree with `standard' CSM classes, because both satisfy the 
Deligne--Grothendieck prescription.

The main new ingredient making our construction possible is the
introduction of {\em `proChow groups',\/} as inverse limits of ordinary
Chow groups over the system of maps to complete varieties. Thus, 
an element of the proChow group $\hA_*U$ is a compatible choice 
of a class in each complete variety to which $U$ maps. The proChow 
group agrees with the conventional Chow group for complete varieties, 
but is in general much larger for noncomplete varieties.

The key feature of proChow groups is that they are functorial with 
respect to arbitrary maps: this is what allows us to define a contribution 
in $A_*X$ from (for example) an {\em open\/} stratum $U$ of $X$. 
If $i_U: U \to X$ is the embedding, there is in general no push-forward 
${i_U}_*: A_*U \to A_*X$, while there {\em is\/} a push forward 
${i_U}_*: \hA_*U\to \hA_*X$ at the level of proChow groups.

Given then the choice of a distinguished element $\bxd U$ in the 
proChow group of every {\em nonsingular\/} variety $U$,
satisfying suitable compatibility properties, we may define an element
$\bxd X\in \hA_*X$ for {\em arbitrary\/} (that is, possibly singular)
varieties by setting
$$\bxd X = \sum_U {i_U}_* \bxd U$$
for any decomposition $X=\amalg_U U$ of $X$ into disjoint nonsingular
subvarieties $U$. If $X$ is complete, this gives a distinguished element 
of the ordinary Chow group $A_*X$ of~$X$.

Describing this mechanism gluing {\em local\/} intersection-theoretic 
information into {\em global\/} one is the second main goal of this article. 
We show (Proposition~\ref{gldprop}) that the compatibility required for 
this definition reduces to a simple blow-up formula. We then prove
(Proposition~\ref{gldcheck}) that this blow-up formula is satisfied by 
the Chern class of the bundle of differential forms with logarithmic poles 
along a divisor at infinity. By the mechanism described above we get a 
distinguished element $\bxd X\in \hA_*X$ for any $X$, and this is our 
{\em pro\/}CSM class of a (possibly singular) variety.

We should point out that the `good local data' arising from the bundle
of differential forms with logarithmic poles is in fact the only
nontrivial case we know satisfying the compatibility requirement of
Proposition~\ref{gldprop}. It would be quite interesting to produce other 
such `gluable' data; perhaps it would be even more interesting to prove
that this is in fact essentially {\em the only\/} such example, as it would 
provide a further sense in which (pro)CSM classes are truly canonical.

The definition of proCSM classes extends immediately to 
{\em constructible functions\/} on $X$, yielding a transformation from 
the functor $F$ of constructible functions to the proChow functor~$\hA_*$. 
Both functors have a push-forward defined for arbitrary (regular) 
morphisms, and we prove (Theorem~\ref{natthm}) that the transformation 
is natural with respect to them. That is: denoting by $\bxd{\varphi}$ the 
proCSM class of the constructible function $\varphi$, we prove that, for 
an {\em arbitrary\/} morphism of varieties $f: X \to Y$,
$$\bxd{f_*(\varphi)} = f_*\bxd{\varphi}\quad.$$
If in particular $X$ and $Y$ are complete, and $f$ is proper, all reduces 
to the more conventional naturality statement and yields a new proof of 
(the Chow flavor of) MacPherson's theorem.

From a technical standpoint, our construction relies on factorization 
of birational maps (\cite{MR2003c:14016}); the fact that this powerful 
result is relatively recent is the likely reason why the construction 
presented here was not proposed a long time ago. Also, we rely on 
MacPherson's `graph construction' in the proof of the key 
Lemma~\ref{keycov}, similarly to MacPherson's own proof of 
naturality in \cite{MR50:13587}.

The relation between CSM classes and Chern classes of bundles
of differential forms with logarithmic poles is not new: 
cf.~Proposition~15.3 in \cite{MR1893006} and Theorem~1
in \cite{MR2001d:14008}. In fact, this relation and MacPherson's
theory may be used to shortcut the paper substantially.
Indeed, our definition of `good local data' may be recast as follows:
for every nonsingular variety $U$, one may define the element
$\bxd U\in \hA_*U$ by selecting, for each complete variety $X$ 
containing $U$, the element
$$c_*(\one_U)\in A_*X$$
obtained by applying MacPherson's natural transformation to the
function that is $1$ over $U$, and $0$ outside of $U$. Both
the basic compatibility (Proposition~\ref{gldcheck}) and the
key Lemma~5.3 follow then from MacPherson's naturality 
theorem. Granting these two facts, the material in \S5 upgrades 
MacPherson's natural transformation $F \to A_*$ (which is natural 
with respect to {\em proper\/} morphisms) to a transformation 
$F \to \hA_*$, natural with respect to arbitrary morphisms.
More details on this approach may be found in \cite{CCSVR}, including 
simple applications (such as a proof of Ehlers' formula 
for the Chern-Schwartz-MacPherson class of a toric variety).

However, working independently of MacPherson's naturality theorem, 
as we do in this paper, allows 
us to discriminate carefully between parts of the construction which 
may extend in a straightforward way to a more general context, and
parts which depend more crucially on (for instance) the characteristic
of the ground field. For example, Proposition~\ref{gldcheck} turns
out to be a purely formal computation, while Lemma~5.3 is much
subtler, and in fact fails in positive characteristic---this distinction 
is lost if one chooses to take the shortcut sketched above. 

Our construction depends on canonical resolution of singularities; 
at the time of this writing, this is only known to hold in characteristic 
zero. Should resolution of singularities be proved in a more general 
setting, our construction will extend to that setting. However, 
characteristic zero is employed more substantially (for example 
through generic smoothness) in the proof of naturality, and simple 
examples show that naturality can{\em not\/} be expected to 
hold in general in positive characteristic. In fact (as J\"org 
Sch\"urmann pointed out to me), covariance of the push-forward for 
constructible functions (Theorem~\ref{covthm}) already fails in positive 
characteristic.

We find this state of affairs intriguing. The construction of 
(pro)CSM classes may well extend to positive characteristic, retaining its 
basic normalization and additivity properties; these `only' depend on 
resolution of singularities, as is shown in this note. But the subtler 
naturality property of these classes cannot carry over, at least within 
the current understanding of the situation.

It is formally possible to extend MacPherson's construction of CSM 
classes to arbitrary characteristic, independently of resolution of 
singularities; for example, this is done in \cite{MR700743}, \S2.5. 
It would be interesting to establish whether proCSM classes agree 
with those defined by Navarro~Aznar. In the presence of naturality 
it is easy to see that proCSM classes agree (for complete varieties, 
and in the Chow group with $\Qbb$ coefficients) with the classes 
discussed in \cite{MR2183846}, \S5; but these latter also 
depend on resolution of singularities.

Similar comments apply to a class which may be defined in general
for hypersurfaces as a twist of Fulton's Chern class, and is known
to agree with the CSM class in characteristic zero (see 
\cite{MR2001i:14009}; both resolution of singularities and naturality
are needed in the proof).

We believe the local-to-global formalism described in this note (maybe
with different target functors rather than Chow) should have other 
applications, for example simplifying the treatment of other 
invariants of singular varieties; this will be explored elsewhere. Clearly
{\em pro-\/}flavors of other functors may be constructed similarly to
the proChow functor presented here; we chose to concentrate on this
example in view of the immediate application to CSM classes. 
J\"org Sch\"urmann has pointed out that it would be worth analyzing the 
construction studied here vis-a-vis the {\em relative Grothendieck group 
of varieties\/} as used in \cite{math.AG/0503492} (particularly in view of 
the parallel between the criterion for `good local data', given here in 
Proposition~\ref{gldprop}, and Franziska Bittner's description of the 
relations defining the Grothendieck group, cf.~\cite{MR2059227}). 
Also, Jean-Paul Brasselet has suggested that the construction of 
proCSM classes may provide an alternative proof of the equality of 
Schwartz and MacPherson classes, that is, the result of 
\cite{MR83h:32011}.

A substantial part of this work was done while the author was visiting the 
{\em Institut de Math\'ematiques de Luminy\/} and the {\em Max-Planck-Institut
f\"ur Mathematik\/} in Bonn, in the summer of 2005. I would especially
like to thank Jean-Paul Brasselet and Matilde Marcolli, for the hospitality 
and for stimulating discussions. I also thank J\"org Sch\"urmann for 
pointing out inaccuracies in an earlier version of this note, and for 
several clarifying remarks concerning the case of positive characteristic;
and Pierre Deligne for comments on an earlier version of this paper.

\newpage


\section{proChow groups}\label{proCg}
 
\subsection{}
We work over an algebraically closed field $k$; further restrictions
will come into play later (\S\ref{reduxdnc} and ff., \S\ref{sumnat} and ff.). 
{\em Schemes\/} will be understood to be separated, of finite type over 
$k$. We say that $X$ is {\em complete\/} if it is proper over $k$.

\subsection{}
We denote by $A_*X$ the conventional Chow group of $X$, as defined 
in \cite{MR85k:14004}; $A_*$ is then a functor from the category of 
schemes to abelian groups, covariant with respect to {\em proper\/}
maps. We begin by defining a functor $\hA_*$ from schemes to 
abelian groups, covariant with respect to arbitrary (regular) morphisms.

For any scheme $U$ consider the category $\mathcal U$ of maps
$$i: U \rightarrow X^i$$
with $X^i$ complete, and morphisms $i\to j$ given by commutative 
diagrams
$$\xymatrix@R=0pt{
& X^i \ar[dd]^\pi \\
U \ar[ur]^i \ar[dr]_j \\
& X^j
}$$
with $\pi$ a proper morphism. 

Note that any two objects $i$, $j$ of $\mathcal U$ are preceded by a third
(that is, $i\times j$), and any two parallel morphisms 
$\xymatrix@C=10pt@1{
i \ar@<2pt>[r] \ar@<-2pt>[r] & j}$
are equalized in $\mathcal U$ (by $U \to \overline{\text{im}(i)}$); that is, 
$\mathcal U$ is a filtered category.

We will say that an embedding $i:U \hookrightarrow X^i$ is a {\em closure\/} 
of $U$ if $X^i$ is complete and $U$ is a dense open set of $X^i$; recall that 
every scheme $U$ has closures (\cite{MR0158892}).

\begin{lemma}\label{reduxlem}
Closures form a small cofinal subcategory $\overline{\mathcal U}$ of 
$\mathcal U$.
\end{lemma}

\begin{proof}
Let $i: U \to X^i$ be any map from $U$ to a complete variety. Let 
$j: U \hookrightarrow \overline U$ be any fixed closure of $U$. Then 
$i$ extends to a {\em rational\/} map 
$\overline i:\xymatrix@1{\overline U \ar@{-->}[r] &  X^i}$, resolving the 
indeterminacies of which produces a diagram
$$\xymatrix{
& \hat U \ar[d] \ar[dr]^{\hat i} \\
U \ar@{^{(}->}[ur]^{\hat j} \ar@{^{(}->}[r]^j \ar@/_1pc/[rr]_i & 
\overline U \ar@{-->}[r]^{\overline i} & X^i
}$$
with $\hat j$ also a closure of $U$: since $i$ is defined on the whole
of $U$, resolving its indeterminacies may be achieved by blowing up
a subscheme in the complement of $U$, so the inclusion $U\subset 
\overline U$ lifts to an inclusion $U\subset \hat U$. The proper map 
$\hat i$ gives a morphism from the inclusion $\hat j$ to $i$, proving 
that closures are cofinal in $\mathcal U$.

The category of closures is small: indeed, by the same reasoning,
any closure of $U$ may be obtained by blowing up a subscheme of 
$\overline U$, then blowing down a subscheme of the blow-up.
\end{proof}

Applying the (conventional) Chow functor to closures of $U$
gives an inverse system of groups 
$\{A_*X^i\}_{i\in Ob(\overline{\mathcal U})}$ 
organized by (proper) push-forward.

\begin{defin}\label{defproC}
The {\em proChow\/} group of $U$ is the inverse limit of this system:
$$\hA_*U:=\varprojlim_{i\in Ob(\overline{\mathcal U})} A_*X^i\quad.$$
\end{defin}

By Lemma~\ref{reduxlem}, we may think of an element $\alpha\in \hA_*U$ 
as the choice of a class $\alpha^i$ in the Chow group of every complete 
variety $X^i$ to which $U$ maps, compatibly with proper push-forward.

\begin{remark}\label{degree}
In particular, elements $\alpha\in\hA_*U$ have a well-defined 
{\em degree\/}
$$\int \alpha \in \Zbb\quad:$$
the structure map $U \to \Spec k$ is a map to a complete variety, so
$\alpha$ determines an element $\int \alpha\in A_*\Spec k=\Zbb$.
\end{remark}

However, the reader should keep in mind that in order to define an element 
of $\hA_*U$ it suffices to define compatible classes for {\em closures\/} of $U$.

\subsection{}\label{pushf}
Any morphism $f: U \to V$ realizes the category $\mathcal V$ 
corresponding to $V$ as a subcategory of $\mathcal U$; thus, a 
compatible assignment of classes in $A_*X^i$ for all $i$ in $\mathcal U$ 
determines in particular a compatible assignment for $i$ in $\mathcal V$. 
That is, $f$ induces a homomorphism
$$f_*: \hA_*U \to \hA_*V\quad.$$
The following remarks should be clear.

\begin{lemma}
With notation as above:
\begin{itemize}
\item $(f\circ g)_*=f_*\circ g_*$: that is, $\hA_*$ is a covariant functor
from the category of algebraic varieties, with morphisms, to
abelian groups.
\item If $X$ is complete, then there is a canonical isomorphism
$\hA_*X\cong A_*X$.
\item If $f: X \to Y$ is a proper map of complete varieties, the
induced homomorphism
$$f_*: A_*X \cong \hA_*X \to \hA_*Y \cong A_*Y$$
is the conventional proper push-forward for Chow groups.
\end{itemize}
\end{lemma}

\subsection{}\label{reduxdnc}
{\em We now assume that canonical resolution of singularities holds.\/}

If $U$ is a {\em nonsingular variety,\/} we say that an embedding 
$i:U \hookrightarrow X^i$ is a {\em good closure\/} of $U$ if it is a 
closure, and further
\begin{itemize}
\item $X^i$ is nonsingular;
\item the complement of $U$ in $X^i$ is a divisor with simple normal 
crossings (that is, normal crossings and nonsingular components).
\end{itemize}

Resolution of singularities implies immediately that good closures
are cofinal in $\overline{\mathcal U}$.

\begin{prop}\label{dnccor}
Assume $U$ is nonsingular.
\begin{itemize}
\item In order to define an element $\alpha\in \hA_*U$ it suffices to assign
$\alpha^i\in A_*X^i$ for good closures $i:U\hookrightarrow X^i$, satisfying
the following compatibility requirement: if $i:U\hookrightarrow X^i$,
$j:U\hookrightarrow X^j$ are good closures, and $i\to j$:
$$\xymatrix@R=0pt{
& X^i \ar[dd]^\pi \\
U \ar[ur]^i \ar[dr]_j \\
& X^j
}$$
with $\pi$ a {\em blow-up of $X^j$ along a smooth center meeting
$X^j\smallsetminus U$ with normal crossings\footnote{As in 
\cite{MR2003c:14016}, this means that the center of the blow-up and
the components of $X^j\smallsetminus U$ are defined (analytically) at each 
point by subsets of a system of local parameters.},\/} then $\alpha^j=\pi_*(\alpha^i)$.
\item Two  elements $\alpha$, $\beta$ in $\hA_*U$ are equal if and only
if $\alpha^i=\beta^i$ for all good closures~$i$.
\end{itemize}
\end{prop}

\begin{proof}
For the first assertion, cofinality implies that it suffices to assign
$\alpha^i$ compatibly for good closures $i$; the reduction to the
case of a blow-up is a standard application of the factorization theorem for 
birational maps (\cite{MR2003c:14016}).

The second assertion is immediate from the cofinality of good closures.
\end{proof}

\subsection{}
A number of further remarks on these notions appear useful in principle,
but will not be needed in this paper.

For example: for every scheme $U$, the surjection from algebraic
cycles $Z_*U$ to the conventional Chow group factors through the
proChow group,
$$\xymatrix@R=10pt@C=0pt{
Z_*U \ar@{->>}[rr] \ar[dr] & & A_*U \\
& \hA_*U \ar@{->>}[ur]
}$$
Indeed, we may associate with every subvariety $W$ of $U$ the class
$[\overline{i(W)}]\in A_*X^i$, for each closure $i$; this
assignment is clearly compatible. In particular, $U$ itself
has a `fundamental class' $[U]\in \hA_*U$.

The surjection $\hA_* U \twoheadrightarrow A_*U$ is compatible with proper 
push-forwards; it may also be realized by mapping $\hA_*U$ to $A_*X^i$ 
for any closure $i$, and following with the natural surjection $A_*X^i 
\twoheadrightarrow A_*U$ (\cite{MR85k:14004},~\S1.8).

If $U$ is nonsingular and $i: U \hookrightarrow X^i$ is a {\em good\/} 
closure, then $\hA_*U \to A_*X^i$ is in fact already a surjection. Indeed,
closures of $U$ dominating $i$ are cofinal in $\mathcal U$; if $\alpha\in A_*X^i$,
assigning to any such closure the pull-back of $\alpha$ defines an element
of $\hA_*U$ mapping to $A_*X^i$. For example, this implies that $\hA_*
\Abb^2$ is not finitely generated. In this sense $\hA_*U$ is, in general, 
`much larger' than $A_*U$.

Variations on the definition of the proChow group arise naturally in certain
contexts: for example, if $U$ admits a torus action then one may take the 
limit of Chow groups as in Definition~\ref{defproC}, but ranging over
{\em toric\/} closures of $U$. Such variations will carry `proCSM classes'
by the same technique introduced in the rest of this paper for the
proChow groups of Definition~\ref{defproC}, provided that certain minimal
requirements are satisfied. Details are left to the interested reader (and 
a slightly more extended discussion may be found in \cite{CCSVR}).


\section{Globalizing local data}\label{gloloc}
 
\subsection{}\label{goodlocdata}
Next we come to the question of defining global invariants on a variety $X$ from
local data: for example, from data given on (open) strata $U$ of a stratification
of $X$. The functorial proChow group offers a natural way to do this:
\begin{itemize}
\item Suppose a class $\bxd U\in \hA_*U$ is defined for every 
{\em nonsingular\/} irreducible variety~$U$;
\item then we will define $\bxd X\in \hA_*X$ by
$$\bxd{X}:=\sum_U {i_U}_* \bxd{U}\quad,$$
where $X$ is the disjoint union of the varieties $U$, each $U$ is nonsingular, 
irreducible, and $i_U:U \to X$ denotes the inclusion.
\end{itemize}
Of course one has to check that this operation is well-defined. We say that
the assignment $U\mapsto \bxd U$ for $U$ nonsingular is {\em good
local data\/} if this is the case. In this section we identify a condition yielding
good local data.

\subsection{}\label{gldcond}
We denote by
$$U\mapsto \fc_U^\oU\in A_*\oU$$
the assignment of a class to each nonsingular variety $U$, in a good closure
$\oU$ of $U$. In our application in \S\ref{proCSM}, $\fc_U^\oU$ will be obtained
as the Chern class of a suitable bundle; in this section we are simply interested
in whether this assignment defines {\em good\/} local data. We will use the 
following notations:
\begin{itemize}
\item $U$: a nonsingular irreducible variety;
\item $\oU$: a good closure of $U$;
\item $D=\oU\smallsetminus U$, a divisor with simple normal crossings in~$\oU$;
\item $W$: a nonsingular closed irreducible subvariety of $\oU$ meeting $D$ 
with normal crossings;
\item $Z$: the intersection $W\cap U$; note that $W$ is a good closure of $Z$
(if $Z\neq\emptyset$);
\item $\pi: \oV\to \oU$: the blow-up of $\oU$ along $W$;
\item $F$: the exceptional divisor $\pi^{-1}(W)$;
\item $\vvv$: the blow-up of $U$ along $Z$, that is, $\pi^{-1}(U)$;
\item $E=\pi^{-1}(Z)=F\cap \vvv$; $F$ is a good closure of $E$ if $Z\neq
\emptyset$.
\end{itemize}
These may be collected in the diagram:
$$\xymatrix@!0{
& \vvv \ar@{^(->}[rr] \ar '[d][dd] & & \oV \ar[dd]^\pi \\
E \ar@{^(->}[ru] \ar@{^(->}[rr] \ar[dd] & & F \ar@{^(->}[ru] \ar[dd]\\
& U \ar@{^(->}'[r][rr] & & \oU\\
Z \ar@{^(->}[ru] \ar@{^(->}[rr] & & W \ar@{^(->}[ru]_w
}$$

\begin{prop}\label{gldprop}
Assume that, for all choices of $U$, $\oU$, etc.~as above, 
$$\fc_U^\oU = \pi_* \fc_{\vvv\smallsetminus E}^{\oV} + w_* \fc_Z^W\quad.$$
Then the assignments $\fc_U^\oU$ determine good local data $\bxd U
\in \hA_*U$ in the sense of \S\ref{goodlocdata}.
\end{prop}

For notational convenience we understand $\fc_Z^W=0$ if $Z=\emptyset$. 
Also note that $\oU$ is complete and $w$, $\pi$ are proper maps; the 
push-forwards appearing in the statement are the ordinary proper 
push-forwards for Chow groups $w_*: A_*W \to A_*\oU$, 
$\pi_*: A_*\oV \to A_*X$. The assignments $U \to \fc_U^\oU$ determine a 
class $\bxd U$ in the proChow group of $U$ by 
Proposition~\ref{dnccor}: it is part of the statement of 
Proposition~\ref{gldprop} that the necessary compatibility is satisfied.

\begin{proof}
If $W$ is disjoint from $U$ then $Z=E=\emptyset$, and $\vvv\cong U$; the 
assumption reduces to
$$\fc_U^\oU = \pi_* \fc_U^{\oV}\quad,$$
that is, the compatibility requirement in Proposition~\ref{dnccor}. Thus the
prescription $U \mapsto \fc_U^\oU$ does define an element $\bxd U$ in
$\hA_*U$.

For $Z\neq\emptyset$, the stated assumption implies that
\begin{equation*}
\tag{*} \bxd U = z_* \bxd Z + i_* \bxd{U\smallsetminus Z}
\end{equation*}
for all nonsingular varieties $U$ and nonsingular closed subvarieties $Z$,
where $z:Z \hookrightarrow U$ and $i:U\smallsetminus Z \hookrightarrow U$
are the embeddings and $z_*$, $i_*$ are the push-forwards defined in
\S\ref{pushf}. Indeed, by Proposition~\ref{dnccor} it suffices to check this
equality after specializing to any good closure of $U$; and we may
dominate any good closure of $U$ with a good closure $\oU$ such that
$W=\overline Z$ is a good closure of $Z$. This is the situation considered
above; since $U\smallsetminus Z\cong \vvv\smallsetminus E$, and $\oV$ is a
good closure of the latter, the formula in the statement of the proposition
implies the claimed equality.

To prove that the assignment $U \mapsto \bxd{U}$ defines good local
data, we have to show that if $X$ is any variety, and 
$$X=\amalg_\alpha \, U_\alpha$$
is a decomposition of $X$ as a finite disjoint union of nonsingular 
subvarieties $\xymatrix@1{U_\alpha \ar@{^(->}[r]^{i_{U_\alpha}} & X}$, 
then
$$\sum_\alpha {i_{U_\alpha}}_* \bxd {U_\alpha}$$
is independent of the decomposition. Any two such decompositions admit
a common refinement, hence we may assume that every element $U$ 
of one decomposition is a finite disjoint union of elements from the other:
$$U=V_1 \amalg V_2 \amalg \cdots \amalg V_r$$
where $U$ and all $V_j$ are nonsingular; and we may assume $V_j$ is 
closed in $V_1\amalg \cdots \amalg V_j$. The required equality
$$\bxd U={i_1}_*\bxd{V_1} + \cdots + {i_r}_*\bxd{V_r}$$
where $i_j:V_i\hookrightarrow U$ denotes the inclusion, is then an immediate 
consequence of (*), concluding the proof.
\end{proof}

\subsection{}
The definition of $\bxd X$ for arbitrary varieties $X$, extending good local 
data for nonsingular varieties as specified in \S\ref{goodlocdata}, satisfies
`inclusion-exclusion'. More precisely:

\begin{prop}\label{incexc}
If $X=\cup_{j\in J} X_j$ is a finite union of subvarieties, then
$$\bxd X = \sum_{\emptyset\neq I\subset J} (-1)^{|I|+1}
{\nu_I}_*\bxd{X_I}$$
where $X_I=\cap_{i\in I}X_i$, and $\nu_I: X_I \hookrightarrow X$ is the
inclusion.
\end{prop}

\begin{proof}
This follows immediately from the case in which $|J|=2$, that is (omitting
push-forwards):
$$\bxd{X\cup Y} = \bxd{X}+\bxd{Y}-\bxd{X\cap Y}\quad,$$
which is straightforward from the definition.
\end{proof}

The same applies {\em a fortiori\/} for the degrees $\int \bxd{X}$
(as defined in Remark~\ref{degree}).

\subsection{}\label{gldcf}
The group of {\em constructible functions\/} of a variety $X$ is the group
$\cF(X)$ of finite integer linear combinations of functions $\one_Z$,
for $Z$ subvarieties of $X$, where
$$\one_Z(p) = \left\{\aligned
1\quad & p\in Z\\
0\quad & p\not\in Z
\endaligned\right.\quad.$$
For $\varphi\in \cF(X)$, we may define an element
$$\bxd{\varphi}\in \hA_*X$$
as follows: if $\varphi=\sum_Z n_Z \one_Z$, we set
$$\bxd{\varphi} = \sum_Z n_Z {\nu_Z}_* \bxd{Z}\quad,$$
where $\nu_Z:Z \hookrightarrow X$ is the inclusion. Inclusion-exclusion
implies that this definition is independent of the decomposition of $\varphi$
(if $\bxd{\cdot}$ arises from good local data).

In good situations, the choice of good local data determines a 
{\em push-forward\/} for constructible functions, as follows: if $f:X \to Y$ is a 
morphism of algebraic varieties, define
$$f_*:F(X) \to F(Y)$$
by setting, for $\varphi=\sum_Z n_Z \one_Z$ as above, and $p\in Y$,
$$f_*(\varphi)(p) = \sum_Z n_Z \int \bxd{f^{-1}(p)\cap Z}\quad;$$
again, the independence on the choices follows from inclusion-exclusion.

For the local data we will define in the next section, and in characteristic zero,
the right-hand-side in this prescription is easily seen to be constructible as
needed (by Lemma~\ref{decom}); and we will show that the resulting 
push-forward is covariant with respect to regular maps. The assignment 
$\varphi \mapsto \bxd{\varphi}$ will then give a transformation of functors
$$\cF \leadsto \hA_*\quad,$$
and the main result in the rest of the paper will be that, for that choice of
local data, this is a {\em natural\/} transformation. 

Specializing to complete varieties, ordinary Chow groups, and proper maps, 
will recover the standard algebraic version of MacPherson's natural 
transformation. 


\section{proCSM classes}\label{proCSM}
 
\subsection{}\label{proCSMld}
We next choose specific local data, by means of Proposition~\ref{gldprop}.

Let $U$ be a nonsingular variety, and let $i:U\hookrightarrow \oU$ be a good 
closure of $U$. In particular, $\oU\smallsetminus U$ is a divisor $D$ with 
normal crossings and nonsingular components $D_i$, $i=1,\dots,r$, in $\oU$.

\begin{defin}
Set
$$\fc_U^\oU:= c(\Omega^1_{\oU}(\log D)^\vee)\cap [\oU]\in A_*\oU\quad.$$
\end{defin}

Here $\Omega^1_{\oU}(\log D)$ denotes the bundle of differential 1-forms
with logarithmic poles along $D$.

As we will prove, this assignment specifies {\em good\/} local data. 
We will use the following immediate lemma:

\begin{lemma}\label{immlem}
Let $\rho^\circ: E \to Z$ be a proper, smooth, surjective map of nonsingular 
varieties. Let $F$, $W$ resp.~be good closure of $E$, $W$, and assume 
that $\rho^\circ$ is the restriction of a proper, smooth, surjective map $\rho: 
F \to W$.
$$\xymatrix{
E \ar[d]_{\rho^\circ} \ar@{^(->}[r] & F \ar[d]^\rho \\
Z \ar@{^(->}[r] & W
}$$
Then
$$\rho_* \fc_E^F = \chi \cdot \fc_Z^W\quad,$$
where $\chi$ is the degree of the top Chern class of the tangent bundle
to any fiber of $\rho$.
\end{lemma}

\begin{proof}
Let $D=W\smallsetminus Z$, a divisor with normal crossings and nonsingular
components by assumption; $F\smallsetminus E=\rho^{-1}(D)$ is then also a 
divisor with simple normal crossings.

In this situation, the exact sequence of differentials on $F$:
$$\xymatrix{
0 \ar[r] &
\rho^*\Omega^1_W \ar[r] &
\Omega^1_F \ar[r] &
\Omega^1_{F|W} \ar[r] & 0
}$$
induces an exact sequence
$$\xymatrix{
0 \ar[r] &
\rho^*\Omega^1_W(\log D) \ar[r] &
\Omega^1_F(\log \rho^{-1}(D)) \ar[r] &
\Omega^1_{F|W} \ar[r] & 0
}$$
The statement follows immediately from this sequence and the projection
formula, since
$$\rho_*\left(c({\Omega^1_{F|W}}^\vee)\cap [F]\right) = \chi\cdot [W]\quad.$$
\end{proof}

\subsection{}
Verifying that the assignment specified above defines good local data is now
a straightforward (but somewhat involved) computation.

\begin{prop}\label{gldcheck}
The classes $\fc_U^\oU$ satisfy the conditions specified in 
Proposition~\ref{gldprop}; therefore, they define good local data in the sense 
of \S\ref{goodlocdata}.
\end{prop}

\begin{proof}
We adopt the notation in \S\ref{gldcond}, and in particular the blow-up
diagram
$$\xymatrix@1{
F \ar[r]^j \ar[d]_\rho & \oV \ar[d]^\pi \\
W \ar[r]_w & \oU
}$$
and we have to verify that
$$\fc_U^\oU = \pi_* \fc_{\vvv\smallsetminus E}^{\oV} + w_* \fc_Z^W\quad.$$
Also recall that
$$c(\Omega^1_{\oU}(\log D)^\vee) = \frac{c(T\oU)}{\prod_i (1+D_i)}$$
(as follows from a residue exact sequence, cf.~\cite{MR98d:32038}, 3.1).

First assume that $Z=\emptyset$, that is, $W$ is contained in $D$.
Denote by $\Til D_i$ the proper transforms of the components $D_i$. 
Then it is easily checked (cf.~Lemma~2.4 in \cite{MR2098642}) that the 
exceptional divisor $F$ and the hypersurfaces $\Til D_i$ together 
form a divisor with simple normal crossings, and their union is the 
complement of $\vvv\cong U$ in $\oV$. The needed statement then 
becomes
$$\pi_*\frac{c(T\oV)}{(1+F)\prod_i (1+\Til D_i)}
\cap [\oV]= \frac{c(T\oU)}{\prod_i (1+D_i)}\cap [\oU]\quad,$$
under the assumption that $W$ is contained in at least one component
of $D$. This is Lemma 3.8, part (5), in \cite{MR2183846}.

If $Z\neq\emptyset$, that is, $W$ is {\em not\/} contained in any
component of $D$, then the proper transforms $\Til D_i$ of $D_i$ agree
with the inverse images $\pi^{-1}(D_i)$, for all $i$. The blow-up $\vvv$ of 
$U$ along $Z$ admits $\oV$ as good closure, with complement 
$\oV\smallsetminus \vvv=\cup \Til D_i$. The projection formula gives
$$\pi_*\left(\frac{c(T\oV)}{\prod_i (1+\Til D_i)}\cap [\oV]\right) 
= \frac 1{\prod_i (1+D_i)}\cap \pi_*(c(T\oV)\cap [\oV])$$
since the class of $\Til D_i$ is $\pi^*(D_i)$. By part (2) of Lemma~3.8
in \cite{MR2183846}, this equals
$$\frac 1{\prod_i (1+D_i)}\cap (c(T\oU)\cap [\oU]+(d-1)
w_* c(TW)\cap [\overline W])\quad,$$
where $d$ denotes the codimension of $W$ in $\oU$ and 
$w:W\hookrightarrow \oU$ is the embedding.

Since $W$ meets $D$ with normal crossings, the divisors $D_i$ cut
out nonsingular divisors on $W$, meeting with normal crossings in $W$, 
and whose union is the complement of $Z$ in $W$. In other words
$W$ is a good closure of $Z$, and the computation given above yields
$$\pi_*\fc_{\vvv}^\oV = \fc_U^\oU + (d-1) w_* \fc_Z^W\quad.$$

On the other hand, $\oV$ is a good closure of $\vvv\smallsetminus E$,
with complement the normal crossing divisor consisting of the
exceptional divisor $F$ and the components $\Til D_i$. Thus
\begin{align*}
\fc_{\vvv\smallsetminus E}^\oV &= 
\frac{c(T\oV)}{(1+F)\prod_i (1+\Til D_i)}\cap [\oV]
=\frac{c(T\oV)}{\prod_i (1+\Til D_i)}\cap [\oV]
-\frac{c(T\oV)}{\prod_i (1+\Til D_i)}\cap [F]\\
&=\frac{c(T\oV)}{\prod_i (1+\Til D_i)}\cap [\oV]
-j_*\frac{c(TF)}{\prod_i (1+j^*\Til D_i)}\cap [F]\quad;
\end{align*}
since $F$ is a good closure of $E$, with complement given by the 
union of the intersections $\Til D_i\cap F$ (of class $j^*\Til D_i$), 
this shows
$$\fc_{\vvv\smallsetminus E}^\oV = \fc_{\vvv}^\oV - j_*\fc_E^F\quad.$$
Combining with the formula obtained above, we get
$$\fc_U^\oU = \pi_*\fc_{\vvv\smallsetminus E}^\oV+w_*\left(\rho_* \fc_E^F -
(d-1)\fc_Z^W\right)\quad.$$
Now $\rho: F \to W$ is a projective bundle, hence smooth and proper,
with fibers $\Pbb^{d-1}$; it restricts to the projective bundle $E \to Z$. 
Applying Lemma~\ref{immlem}, with $\chi=\int c(T\Pbb^{d-1})\cap 
[\Pbb^{d-1}]=d$, gives
$$\rho_*\fc_E^F=d\, \fc_Z^W\quad,$$
concluding the proof.
\end{proof}

\subsection{}
We are ready to define {\em proCSM classes,\/} and the corresponding
transformation $F \leadsto \hA_*$.

\begin{defin}\label{proCSMdef}
The {\em proCSM class\/} of a (possibly singular) variety $X$ is the
class 
$$\bxd{X}\in \hA_*X$$
in the proChow group of $X$, defined by patching the local data
defined in \S\ref{proCSMld}, as explained in \S\ref{gloloc}.
\end{defin}

Explicitly, write $X=\amalg_\alpha U_\alpha$ in any way as  a finite
disjoint union of nonsingular subvarieties $\xymatrix@1{U_\alpha 
\ar@{^(->}[r]^{i_{U_\alpha}} & X}$; then $\bxd X:=\sum_\alpha
{i_{U_\alpha}}_* \bxd{U_\alpha}$ is independent of the 
decomposition, by Proposition~\ref{gldcheck}.

The following `normalization' properties are easy consequences of the
definition.
\begin{prop}\label{norm}
\begin{itemize}
\item If $X$ is complete and nonsingular, then
$$\bxd{X}=c(TX)\cap [X]\quad.$$
\item If $X$ is a compact complex algebraic variety, then
$$\int{\bxd{X}}=\chi_{\text{top}}(X)\quad,$$
the topological Euler characteristic of $X$.
\end{itemize}
\end{prop}

\begin{proof}
The first statement is immediate, as $X$ is a good closure of itself and 
$\hA_*X=A_*X$ if $X$ is complete, and $\fc_X^X=c(TX)\cap [X]$.

For the second statement, since both $\int\bxd{\cdot}$ and 
$\chi_{\text{top}}$ satisfy inclusion-exclusion
we only need to check this equality for $X$ compact {\em and
nonsingular.\/} By the first statement (and the Poincar\'e-Hopf
theorem)
$$\int\bxd{X}=\int c(TX)\cap [X] = \chi_{\text{top}}(X)$$
in this case, as needed.
\end{proof}

\subsection{}
The proCSM class {\em of a constructible function\/} $\varphi$ on a
variety $X$, $\bxd{\varphi}$, and the push-forward of constructible
functions $f_*$ may now be defined as in \S\ref{gldcf}.

The second statement in Proposition~\ref{norm} implies that, for compact 
complex algebraic varieties and $f$ proper, this definition of push-forward
for constructible functions agrees with the conventional one (as given in 
\cite{MR85k:14004}, \S19.1.7).

However, covariance properties of the more general push-forward
introduced here, and the naturality of the corresponding transformation
$F \leadsto \hA_*$, do not appear to be immediate. We will address these
questions in the next section.


\section{The natural transformation $F \leadsto \hA_*$}\label{nattr}

\subsection{}\label{sumnat}
To summarize, we have now defined a notion of {\em proCSM class\/} for
possibly noncomplete, possibly singular varieties $X$, in the proChow group 
of $X$:
$$\bxd{X}\in \hA_*X\quad.$$
This definition relies on the local-to-global machinery of 
\S\ref{gloloc}. We have used resolution of singularities, and the factorization 
theorem of \cite{MR2003c:14016}; the definition of proCSM class can be 
given in any context in which these tools apply.

By contrast, the statement to be proved in this section will use characteristic
zero more crucially. {\em Therefore, in the rest of the paper we work over
an algebraically closed field of characteristic zero.\/}

In \S\ref{gldcf} we have also proposed a notion of push-forward of 
constructible functions
$$f_*:F(X) \to F(Y)$$
for any morphism $f: X \to Y$, and a transformation
$$F \leadsto \hA_*\quad.$$
This depends on $f_*(\one_Z)$ being a constructible function on $Y$, for
any subvariety $Z$ of $X$; this is easily seen to be the case for the definition
obtained from the data defined in \S\ref{proCSM}, at least in characteristic~zero.

\begin{theorem}[Covariance]\label{covthm}
Let $f:X \to Y$, $g:Y \to Z$ be morphisms of algebraic varieties. Then
$$(g\circ f)_* = g_*\circ f_*$$ 
as push-forwards $F(X) \to F(Z)$.
\end{theorem}

\begin{theorem}[Naturality]\label{natthm}
The transformation $F\leadsto \hA_*$ is a {\em natural\/} transformation
of covariant functors from the category of algebraic varieties over
an algebraically closed field of characteristic zero, with morphisms,
to the category of abelian groups.
\end{theorem}

This is the main result of this paper.
In view of Proposition~\ref{norm}, Theorem~\ref{natthm} shows that there 
is a natural transformation from $F$ to $\hA_*$, and hence to homology, 
specializing to the total Chern class of the tangent bundle on nonsingular 
varieties. This recovers Theorem~1 in \cite{MR50:13587}. The uniqueness
of the natural transformation is immediate from resolution of singularities, 
and it follows that the (image in homology of the) proCSM classes defined 
in \S\ref{proCSM} agree with the classes constructed by 
MacPherson. This also follows more directly from MacPherson's
theorem, arguing essentially as in the proof of Proposition~\ref{norm}. 

In order to provide a self-contained treatment of (pro)CSM classes, 
MacPherson's theorem will not be used in the proof of 
Theorems~\ref{covthm} and~\ref{natthm} (MacPherson's
{\em graph construction\/} will be an ingredient in the proof of 
Lemma~\ref{keycov}). 

\subsection{}\label{charp}
It is natural to wonder whether the characteristic~zero restriction in 
Theorems~\ref{covthm} and~\ref{natthm} is crucial, or whether it is a 
technical requirement for our approach to the proof. The restriction
is in fact necessary: as J\"org Sch\"urmann pointed out to me, the
presence in characteristic $p>0$ of \'etale self-covers of 
$\Abb^1$, such as the Artin-Schreier map $x\mapsto x^p-x$, gives a 
counterexample to covariance. Other simple examples (such as the 
Frobenius map) indicate that the difficulty cannot be circumvented by 
naive modifications to the definition of push-forward of constructible 
functions.

We do not know, even at a conjectural level, how the formalism could
be modified in order to avoid the characteristic zero requirement. In our
argument the condition will enter in the proof of Lemma~\ref{keycov}
(specifically, in the key Lemma~\ref{surj}), and through `generic 
gentleness'. Sch\"urmann's example shows that Lemma~\ref{keycov} 
fails in positive characteristic.

\subsection{}
The proofs of Theorems~\ref{covthm} and~\ref{natthm} depend on the
following lemma, which at first sight would appear to be a harmless 
generalization of the trivial Lemma~\ref{immlem}. Its proof is on the
contrary rather technical, and we postpone it to \S\ref{keycovpf}.
We state the lemma here, and use it to prove Theorems~\ref{covthm} 
and~\ref{natthm} in the rest of this section.

\begin{lemma}\label{keycov}
Let $f:U \to V$ be a proper, smooth, surjective map of nonsingular 
varieties over an algebraically closed field of characteristic~zero. Then
$$f_*\bxd{U} = \chi_f\cdot \bxd{V}\quad,$$
where $\chi_f=\int\bxd{f^{-1}(p)}$ (for any $p\in V$).
\end{lemma}

Since by hypothesis the fibers $f^{-1}(p)$ are complete and nonsingular,
$\chi_f$ equals the degree of the top Chern class of their tangent bundle
(by Proposition~\ref{norm}).

\subsection{}
First, we formulate an upgrade of Lemma~\ref{keycov} to a mildly larger
class of maps.

\begin{defin}\label{gentle}
A morphism $f: U \to V$ of nonsingular varieties is  {\em gentle\/} if 
it is smooth and surjective, and further there is a variety $\underline U$ 
and a {\em proper,\/} smooth, surjective morphism 
$\underline f: \underline U \to V$ such that:
\begin{itemize}
\item $U$ is an open dense subset of $\underline U$, and $\underline f|_U
=f$;
\item the complement $H=\underline U\smallsetminus U$ is a divisor with 
normal crossings and nonsingular components $H_i$;
\item letting $H_I$ denote the intersection $\cap_{i\in I} H_i$ (so
$H_\emptyset = \underline U$, and each $H_I$ is nonsingular), each 
restriction
$$\underline f|_{H_I}: H_I \to V$$
is proper, smooth, and surjective.
\end{itemize}
\end{defin}

\begin{lemma}\label{gentle1}
If $f: U\to V$ is gentle, then the number
$$\chi_f:=\int \bxd{f^{-1}(p)}$$
is independent of $p\in V$; further, in characteristic zero, 
$$f_*\bxd U = \chi_f \cdot \bxd V\quad.$$
\end{lemma}

\begin{proof}
By inclusion-exclusion (Proposition~\ref{incexc}):
$$\bxd{U} = \sum_I (-1)^{|I|+1} {i_{H_I}}_*\bxd{H_I}$$
and
$$\int\bxd{f^{-1}(p)} = \sum_I (-1)^{|I|+1} \int
\bxd{H_I\cap \underline f^{-1}(p)}$$
with notation as in Definition~\ref{gentle}, and denoting by $i_{H_I}$ the 
inclusion $H_I\hookrightarrow \underline U$.

Since each $H_I$ maps properly, smoothly, and surjectively onto $V$, the
statements follow from Lemma~\ref{keycov}. 
\end{proof}

\subsection{}
We also single out the following easy properties of `gentleness':
\begin{lemma}\label{gentle2}
\begin{itemize}
\item If $f:U\to V$ is gentle, and $W\subset V$ is a nonsingular subvariety, 
then the restriction
$$f^{-1}(W) \to W$$
is gentle;
\item `Generic gentleness' holds in characteristic zero: if $f: X \to Y$ is any
morphism of varieties, and $X$ is nonsingular, then there 
exists a nonempty, nonsingular open subset $V\subset Y$ such that $f$ 
restricts to a gentle map $U=f^{-1}(V) \to V$;
\item If $f$, $g$, and $g\circ f$ are all gentle:
$$\xymatrix{
U \ar[r]^f \ar@/_1pc/[rr]_{g\circ f} & V \ar[r]^g & W
}$$
then $\chi_{g\circ f}=\chi_g\cdot \chi_f$.
\end{itemize}
\end{lemma}

\begin{proof}
The first statement is clear.

The second is straightforward from Nagata's theorem extending $f$ to a
proper map (\cite{MR0158892}), embedded resolution of singularities
to guarantee the complement is a divisor $H$ with simple normal 
crossings, and ordinary generic {\em smoothness\/} 
(Corollary 10.7 in \cite{MR0463157}) applied to all intersections of 
components of~$H$.

For the third, let $p$ be any point of $W$; by the first statement, the
restriction
$$f_p: (g\circ f)^{-1}(p) \to g^{-1}(p)$$
is gentle (as $f$ is); therefore, by Lemma~\ref{gentle1},
$${f_p}_* \bxd{(g\circ f)^{-1}(p)} = \chi_f\cdot \bxd{g^{-1}(p)}
\quad.$$
Taking degrees gives the stated equality, since $g\circ f$ and $g$ are
gentle, and push-forwards preserve degrees.
\end{proof}

\subsection{}
Generic gentleness lets us decompose {\em any\/} map into gentle ones
(in characteristic zero):

\begin{prop}\label{splitup}
Let $f:X \to Y$ be any morphism of varieties over an algebraically
closed field of characteristic zero. Then there are decompositions 
$$X = \amalg_{\alpha, i} U_{\alpha i} \quad,
\quad Y = \amalg_\alpha V_\alpha$$
into disjoint nonsingular irreducible subvarieties, such that, for all 
$\alpha$ and $i$, $f$ restricts to a gentle map
$$f_{\alpha i} := f|_{U_{\alpha i}}: U_{\alpha i} \to V_\alpha\quad.$$ 
\end{prop}

\begin{proof}
This follows immediately from generic gentleness, after decomposing $X$ 
as a disjoint union of nonsingular irreducible subvarieties.
\end{proof}

\subsection{}
The proofs of Theorem~\ref{covthm} and~\ref{natthm} are now 
straightforward. We begin with covariance.

\begin{proof}[Proof of Theorem~\ref{covthm}]
It suffices to show that the two push-forwards agree on the characteristic
function of each subvariety of $X$, and by restricting $f$ we are reduced
to showing
$$(g\circ f)_*(\one_X)=g_*\circ f_*(\one_X)\quad.$$
Two applications of Proposition~\ref{splitup} yield decompositions
$$X=\amalg_{\alpha, i, j} U_{\alpha i j}\quad, \quad
Y=\amalg_{\alpha,i} V_{\alpha i}\quad,
\quad Z=\amalg_{\alpha} W_\alpha$$
such that all restrictions $U_{\alpha i j} \to V_{\alpha i}$, $V_{\alpha i}
\to W_\alpha$ are gentle. Generic gentleness allows us to assume that
the compositions $U_{\alpha i j} \to W_\alpha$ are also gentle. Denote
by $\chi'_{\alpha i j}$, $\chi_{\alpha i}$, $\chi_{\alpha i j}$ the
corresponding fiberwise degrees; by Lemma~\ref{gentle2}, 
$$\chi_{\alpha i j} =\chi'_{\alpha i j}  \cdot \chi_{\alpha i} \quad.$$

The computation is then completely straightforward:
$$
g_*(f_*(\one_{U_{\alpha i j}}))
= g_* (\chi'_{\alpha i j} \one_{V_{\alpha i}})
= \chi'_{\alpha i j} \cdot \chi_{\alpha i} \one_{W_\alpha}
= \chi_{\alpha i j} \one_{W_\alpha}
=(g\circ f)_* (\one_{U_{\alpha i j}})\quad,
$$
and $g_*(f_*(\one_X))=(g\circ f)_*(\one_X)$ follows by linearity as
$\one_X =\sum_{\alpha, i, j} \one_{U_{\alpha i j}}$.
\end{proof}

\subsection{}
For naturality, note that any splitting of a morphism $f:X \to Y$ into gentle maps 
gives a parallel splitting of both $f_*(\one_X)$ and $f_*\bxd{X}\,$:

\begin{lemma}\label{decom}
Let $f:X \to Y$ be a morphism of varieties in characteristic zero,
and let $X=\amalg U_{\alpha i}$, $Y=\amalg V_\alpha$ be decompositions
as in Proposition~\ref{splitup}. For each~$\alpha$, let 
$\chi_\alpha=\sum_i \chi_{f_{\alpha i}}$.
Then
$$f_*(\one_X) = \sum_{\alpha\in A} \chi_\alpha \one_{V_\alpha}\quad;$$
$$f_*\bxd X = \sum_{\alpha\in A} \chi_\alpha \cdot \bxd{V_\alpha}
\quad.$$
\end{lemma}

\begin{proof}
Any given $p\in Y$ is in precisely one $V_\alpha$; and then $f^{-1}(p)$
is the disjoint union of the fibers of $f_{\alpha i}$. Hence
$$\int \bxd{f^{-1}(p)} = \sum_i \int \bxd{f_{\alpha i}^{-1}(p)}
=\sum_i \chi_{f_{\alpha i}} = \chi_\alpha\quad.$$
This gives the first formula, by definition of push-forward of constructible
functions.

The second formula follows from Lemma~\ref{gentle1}:
$$f_*\bxd X= \sum_{\alpha, i} f_*\bxd{U_{\alpha i}} =
\sum_\alpha \sum_i \chi_{f_{\alpha i}}\cdot \bxd{V_\alpha} =
\sum_\alpha \chi_\alpha \cdot \bxd{V_\alpha}\quad,$$
since each $f_{\alpha i}:U_{\alpha i} \to V_\alpha$ is gentle.
\end{proof}

Naturality follows immediately:

\begin{proof}[Proof of Theorem~\ref{natthm}]
We have to show that for any $f:X \to Y$, and any constructible function
$\varphi\in F(X)$,
$$f_*\bxd\varphi = \bxd{f_*\varphi}\quad.$$
By linearity of (both) $f_*$ it suffices to prove this equality for the 
characteristic function of any subvariety of $X$; by restricting $f$ we may 
assume this subvariety is $X$. That is, it suffices to prove that
$$f_*\bxd {\one_X} = \bxd{ f_*(\one_X)}\quad;$$
and this follows immediately from Lemma~\ref{decom}.
\end{proof}

\subsection{}
This concludes the verification of the Deligne-Grothendieck conjecture
for the local data defined in \S\ref{proCSM}.
The only outstanding item is the proof of the key Lemma~\ref{keycov},
with which we will close the paper.


\section{Proof of Lemma~\ref{keycov}}\label{keycovpf}

\subsection{}
We have to prove that if $U$, $V$ are nonsingular varieties over an
algebraically closed field of characteristic zero, and $f:U \to V$
is proper, smooth, and surjective, then
$$f_*\bxd U = \chi_f \cdot \bxd V\quad,$$
where $\chi_f$ is the degree of the proCSM class of the fibers of~$f$;
that is, $\chi_f$ equals the degree of the top Chern class of the tangent
bundle of any fiber of $f$.

By Proposition~\ref{dnccor} and the definition of push-forward on 
proChow groups, this amounts to verifying that, for any good closure 
$j: V \hookrightarrow Y$,
$${\bxd U}^{j\circ f} = \chi_f\cdot {\bxd V}^j\quad.$$
Dominating $j\circ f$ by a good closure $i:U \hookrightarrow X$, we 
have the following fiber square:
$$\xymatrix{
U \ar[r]^i \ar[d]_f & X \ar[d]^g \\
V \ar[r]_j & Y
}$$
with $i$ and $j$ good closures, $f$ smooth, $f$ and $g$ proper and 
surjective, and we have to verify that
$$g_* \fc_U^X = \chi_f \cdot \fc_V^Y\quad.$$
With the local data chosen in \S\ref{proCSMld}, this amounts to the
following claim:

\begin{claim}\label{keyclaim}
In a fiber square as above, denote by $D$, $E$ the complements 
$X\smallsetminus U$, $Y\smallsetminus V$ (which are divisors with 
simple normal crossings by assumption). Then
$$g_*\left(c(\Omega^1_X(\log D)^\vee)\cap [X]\right)
=\chi_f \cdot c(\Omega^1_Y(\log E)^\vee)\cap [Y]\quad.$$
\end{claim}

This is our objective. The reader should compare Claim~\ref{keyclaim}
with Lemma~\ref{immlem}: the given formula is immediate
if $g$ is a {\em smooth\/} map extending the smooth map $f$; while the
general case of Claim~\ref{keyclaim} requires a bit of work (presented 
in the next several subsections) and relies more substantially on the 
hypothesis on the characteristic.

Claim~\ref{keyclaim} is in a sense equivalent to the naturality
of (pro)CSM classes: Theorem~\ref{natthm} is proved in \S\ref{nattr}
as a consequence of Lemma~\ref{keycov} (and hence of 
Claim~\ref{keyclaim}); conversely, Claim~\ref{keyclaim} could be 
proved as a corollary of MacPherson's naturality theorem (exercise
for the reader!). In order to keep the paper self-contained, the proof 
given here does not assume the result of \cite{MR50:13587}.

\subsection{}\label{summ}
We will prove Claim~\ref{keyclaim} by applying the {\em graph 
construction\/} (see \cite{MR85k:14004}, Chapter~18 and especially 
Example~18.1.6) to the (logarithmic) differential map
$$dg:g^*\Omega^1_Y(\log E) \to \Omega^1_X(\log D)\quad.$$
As $g$ is smooth along $U$, this map is injective over $U$.
The graph construction produces a `cycle at infinity' measuring the
singularities of $dg$, which may be used to evaluate the difference 
in the Chern classes of Claim~\ref{keyclaim}. We will show that all
components of the cycle at infinity dominating loci within 
$D=X\smallsetminus U$ give vanishing contribution to this difference, and 
the claim will follow.

\subsection{}
For notational convenience we switch to bundles of differential
forms with logarithmic poles (rather than their duals). The formula
in Claim~\ref{keyclaim} is equivalent to
$$g_*\left(c(\Omega^1_X(\log D))\cap [X]\right)
=\underline\chi_f \cdot c(\Omega^1_Y(\log E))\cap [Y]\quad,$$
where $\underline\chi_f$ is the degree of the top Chern class of
the {\em co\/}tangent bundle of any fiber of~$f$.

Denote by $m$, $n$ resp.~the dimensions of $X$, $Y$. The 
differential $dg$ determines a {\em rational\/} map to the Grassmann
bundle
$$\xymatrix{
\gamma: X\times \Pbb^1 \ar@{-->}[r] & \Gbb
=\Grass_n(g^*\Omega^1_Y(\log E) \oplus \Omega^1_X(\log D))
}$$
restricting on $X\times \{(\lambda:1)\}$ for $\lambda\ne 0$ to the map
assigning to $x\in X$ the graph of $\frac 1\lambda dg$ at $x$. Over 
$X\times \infty:=X\times \{(1:0)\}$ (that is, for $\lambda\to\infty$) 
$\gamma$ acts as the section corresponding to $g^*\Omega^1_Y(\log E) 
\oplus 0$.

\subsection{}
The indeterminacies of $\gamma$ are contained in $X\times 0:=
X\times \{(0:1)\}$; and in fact in $D\times 0\subset X\times 0$
(since $dg$ is injective along $U$). Closing the graph of $\gamma$
(or equivalently blowing up the ideal of indeterminacies) gives
a variety $\Til{X\times\Pbb^1}$ and a {\em regular\/} lift $\tgamma$
of $\gamma$:
$$\xymatrix@R=15pt{
\Til{X\times \Pbb^1} \ar@/_2pc/[ddd]_p \ar[d]^\pi \ar[dr]^{\tgamma} \\
X\times \Pbb^1 \ar[d]^\rho \ar@{-->}[r]^\gamma & \Gbb \ar[dl] \\
X \ar[d]^g \\
Y
}$$
Now
$$[\pi^{-1}(X\times \infty)] = [\pi^{-1}(X\times 0)]$$
as rational equivalence classes of divisors. Note that $\pi^{-1}(X\times
\infty)$ maps isomorphically to $X$. As for $\pi^{-1}(X\times 0)$,
it consists of the proper transform $\Til X$ of $X\times 0$ and of 
the components $\Gamma_i$ of the exceptional divisor, appearing with 
multiplicities $r_i$.

\subsection{}
Let $\cQ$ denote the {\em universal quotient bundle\/} over $\Gbb$, 
a bundle of rank~$m$. We have
$$c(\tgamma^*\cQ)\cap [\pi^{-1}(X\times \infty)]=
c(\tgamma^*\cQ)\cap \left([\Til X]+ \sum_i r_i 
[\Gamma_i]\right)$$
in $A_*(X\times\Pbb^1)$.

\begin{lemma}\label{twoeq}
The following equalities hold in $Y$:
\begin{itemize}
\item $p_*\left(c(\tgamma^*\cQ)\cap [\pi^{-1}(X\times \infty)]
\right)=g_*\left(c(\Omega^1_X(\log D))\cap [X]\right)\quad$;
\item $p_*\left(c(\tgamma^*\cQ)\cap [\Til X]\right) = 
\underline\chi_f\cdot c(\Omega^1_Y(\log E))\cap [Y]\quad$.
\end{itemize}
\end{lemma}

\begin{proof}
These are easy consequences of the basic set-up. For the second
equality, the restriction $\tgamma'$ of $\tgamma$ to $\Til X$ factors through
$$\Gbb':=\Grass_n(\Omega^1_X(\log D))\cong 
\Grass_n(\Omega^1_X(\log D)\oplus 0)
\subset \Gbb\quad.$$
Over $\Gbb'$, $\cQ$ splits as the direct sum of the universal quotient 
bundle $\cQ'$ of rank $(m-n)$ of $\Gbb'$, and the pull-back of 
$g^*\Omega^1_Y(\log E)$; thus
$$c(\tgamma^*\cQ)\cap [\Til X]=c(p^*\Omega^1_Y(\log E))\cdot
c(\tgamma'^*\cQ')\cap [\Til X]\quad,$$
and the given formula follows immediately, since $\cQ'$ restricts to
the cotangent bundle on fibers over points of $V$.
\end{proof}

\subsection{}
As promised in the short summary in \S\ref{summ}, Lemma~\ref{twoeq} 
shows that the difference between the classes appearing in 
Claim~\ref{keyclaim} is controlled by components of the cycle at infinity 
in the graph construction. Explicitly:

\begin{corol}
$$g_*\left(c(\Omega^1_X(\log D))\cap [X]\right)
-\underline\chi_f \cdot c(\Omega^1_Y(\log E))\cap [Y]
=\sum_i r_i \cdot p_*\left(c(\cQ)\cap [\Gamma_i]\right)$$
\end{corol}

Therefore, the following claim will conclude the proof of 
Lemma~\ref{keycov}:

\begin{claim}\label{gammato0}
For any component $\Gamma$ of  the exceptional divisor in 
$\Til{X\times \Pbb^1}$, 
$$p_*\left(c(\cQ)\cap [\Gamma]\right) = 0\quad.$$
\end{claim}

Incidentally, up to this point the discussion could have been carried out
for the ordinary differential of the map $g$; but Claim~\ref{gammato0}
fails for the ordinary differential (cf.~\cite{MR85k:14004}, 
Example~18.1.6~(f)). Claim~\ref{gammato0} is a remarkable property 
of the {\em logarithmic\/} differential (in characteristic zero).

\subsection{}
The proof of Claim~\ref{gammato0} will rely on the following general
observation.

\begin{lemma}\label{vanlem}
Let $p: \Gamma \to W$ be a proper morphism of schemes, and let $\cQg$ 
be a vector bundle on $\Gamma$, of rank $\le \dim \Gamma$. 
Assume that there is a surjection
$\xymatrix{
\cQg \ar@{->>}[r] & p^*\cT
}$
of coherent sheaves on $\Gamma$, where $\cT$ is a coherent sheaf
on $W$ of rank $>\dim W$. Then
$$p_*\left(c(\cQg)\cap [\Gamma]\right) = 0\quad.$$
\end{lemma}

\begin{proof}
There exists (see for example~\cite{MR700743}, \S1.1) a proper birational
morphism $\nu: \Til W \to W$ such that the pull-back $\nu^*({\mathcal T})$ 
is locally free modulo torsion: in particular, there is a surjection
$\xymatrix@1{
\nu^*\cT \ar@{->>}[r] & \hcT
}$
with $\hcT$ locally free on $\Til W$, of rank equal to the rank of 
$\cT$ (and hence $>\dim W$). Let $\tGamma$ be the component
dominating $\Gamma$ in the fiber product:
$$\xymatrix{
\tGamma \ar[d]_{\widetilde p} \ar[r]^{\widehat\nu} & \Gamma \ar[d]^p \\
\Til W \ar[r]_\nu & W
}$$
Then we have surjections
$$\xymatrix{
\widehat\nu^*\cQg \ar@{->>}[r] & \widehat\nu^* p^*\cT =
\widetilde p^* \nu^*\cT \ar@{->>}[r] & \widetilde p^* \hcT
}$$
with $\widetilde p^*\hcT$ locally free on $\tGamma$, of rank $>\dim W$.
Let $\cK$ be the kernel, so we have the exact sequence of locally free
sheaves on $\tGamma$:
$$\xymatrix{
0 \ar[r] & \cK \ar[r] & \widehat\nu^*\cQg \ar[r] & \widetilde p^* \hcT \ar[r]
& 0
}$$
and $\rk\cK=\rk \cQg-\rk \hcT<\dim\tGamma - \dim\Til W$. It follows that
$\widetilde p_* (c(\cK)\cap [\tGamma]) = 0$, and this implies the
statement by the projection formula: $p_*(c(\cQ)\cap [\Gamma])$ equals
$$p_* \widehat\nu_* (c(\widehat\nu^*\cQ)\cap [\tGamma]) 
= \nu_*\widetilde p_*\left( c(\widetilde p^*\hcT)\cdot c(\cK)\cap[\tGamma]
\right)
= \nu_* c(\hcT)\cap \widetilde p_* \left(c(\cK)\cap [\tGamma]\right) = 0$$
as needed.
\end{proof}

\subsection{}
We choose a component $\Gamma$ of the exceptional divisor, and
restrict all relevant maps to it:
$$\xymatrix{
\Gamma \ar@/_1pc/[dd]_p \ar[d]^\sigma \ar[r]^{\tgamma} & \Gbb 
\ar[d] \\
Z \ar[d] \ar@{^(->}[r] & X \ar[d]^g \\
W \ar@{^(->}[r] & Y
}$$
Here $Z$ and $W$ are the images of $\Gamma$ in $X$, $Y$ respectively.
Note that $Z\subset D$ and $W\subset E$; we assume $W\subset E_i$
for $i\le s$, $W\not\subset E_i$ for $i>s$, and we denote by 
$\Eus$ the intersection $E_1\cap\cdots\cap E_s$.

We let $\cSg$, $\cQg$ be resp.~the pull-back to $\Gamma$ of the universal 
sub-~and quotient bundle over $\Gbb$; hence we have an exact sequence
$$\xymatrix{
0 \ar[r] & \cSg \ar[r] & \sigma^* (g^*\Omega^1_Y(\log E) \oplus 
\Omega^1_X(\log D))|_Z \ar[r] & \cQg \ar[r] & 0
}\quad.$$

\subsection{}
The residue exact sequence for the bundle of logarithmic differential forms
restricts to an exact sequence
$$\xymatrix{
0 \ar[r] & \Omega^1_{\Eus}|_W \ar[r] & \Omega^1_Y(\log E)|_W \ar[r]
& \mathcal O_{W}^{\oplus s}\oplus (\oplus_{i>s} 
\mathcal O_{Z\cap E_i}) \ar[r] & 0
}\quad.$$
In particular, $\Omega^1_Y(\log E)|_W$ contains a distinguished copy of 
the {\em conormal sheaf\/} $N^*_W\Eus$, that is, the kernel of the
natural surjection $\Omega^1_{\Eus}|_W \to \Omega^1_W$.
The quotient defines a coherent sheaf $\cT$ on $W$, which by
construction fits in an exact sequence
$$\xymatrix{
0 \ar[r] & \Omega^1_W \ar[r] & {\mathcal T} \ar[r] & \mathcal O_W^{\oplus s} 
\oplus (\oplus_{i>s} \mathcal O_{Z\cap E_i}) \ar[r] & 0
}$$
In particular, note that $\rk\cT>\dim W$.

\begin{lemma}\label{surj}
There is a surjection
$$\xymatrix{
\cQg \ar@{->>}[r] & p^*\cT
}$$
of coherent sheaves on $\Gamma$.
\end{lemma}

\begin{proof}
It suffices to show that the image of $\cSg$ in 
$$\sigma^*(g^*\Omega^1_Y(\log E)\oplus 0)|_Z\cong 
p^*\Omega^1_Y(\log E)|_W$$ 
is contained in the image of $p^* N^*_W \Eus$:
$$\xymatrix{
0 \ar[r] & \cSg \ar[r] & \sigma^*(g^*\Omega^1_Y(\log E) \oplus 
\Omega^1_X(\log D))|_Z \ar[d] \ar[r] & \cQg \ar[r] \ar@{-->}[d] & 0 \\
& p^*N^*_W \Eus \ar[r] & p^*\Omega^1_Y(\log E)|_W \ar[r] & 
p^*{\mathcal T} \ar[r] & 0
}$$
and this may be verified by a computation in local coordinates.
\end{proof}

\subsection{}
By Lemma~\ref{vanlem}, Lemma~\ref{surj} implies the vanishing prescribed
in Claim~\ref{gammato0}, concluding the proof of Lemma~\ref{keycov}.

The coordinate computation in Lemma~\ref{surj} uses characteristic zero: 
for example, in characteristic~$p>0$ problems arise if some component of 
$D$ dominating a component of $E$ appears with multiplicity equal to a 
multiple of $p$. This is in fact precisely what happens with the Artin-Schreier 
map, which gives (as mentioned in \S\ref{charp}) a simple counterexample 
to Lemma~\ref{keycov} in positive characteristic.



\end{document}